\documentclass{pspum-l}
\usepackage{amscd, amssymb}

\theoremstyle{plain}

\newtheorem{lemma}{Lemma}
\newtheorem{theorem}{Theorem}
\newtheorem*{theorem*}{Theorem}

\newtheorem*{claim*}{Claim}

\newtheorem*{conjecture*}{Conjecture}
\newtheorem*{question*}{Question}

\newtheorem*{Virasoro}{Virasoro Conjecture}

\theoremstyle{definition}
\newtheorem{definition}{Definition}
\newtheorem*{definition*}{Definition}

\theoremstyle{remark}
\newtheorem*{remark*}{Remark}
\newtheorem*{remarks*}{Remarks}

\newtheorem*{example*}{Example}
\newtheorem*{acknowledgments}{Acknowledgments}

\newcommand{\PP}{\mathbb{P}}
\newcommand{\QQ}{\mathbb{Q}}

\newcommand{\Mbar}{{\overline{M}}}

\newcommand{\p}{\partial}

\newcommand{\ev}{{\operatorname{ev}}}
\newcommand{\ft}{\operatorname{ft}}

\newcommand{\vir}{\operatorname{vir}}
\newcommand{\la}{\langle}
\newcommand{\ra}{\rangle}
\newcommand{\on}{\operatorname}

\newcommand{\HH}{\mathcal{H}}

\title[Axiomatic Gromov--Witten theory]
{Notes on Axiomatic Gromov--Witten theory and applications}

\author{Y.-P.~Lee}

\address{Department of Mathematics \\
        University of Utah \\
        Salt Lake City, Utah 84112-0090\\
        U.S.A.}
\email{yplee@math.utah.edu}

\begin{document}


\maketitle

\setcounter{section}{-1}

\section{Introduction} \label{s:1}

The purpose of these notes is to give their readers some idea of Givental's 
\emph{axiomatic Gromov--Witten theory}, and a few applications.
Due to the scope of these notes, some statements are not 
precisely formulated and almost all proofs are omitted.
However, we try to point out some subtleties, and to give references for 
further reading whenever desirable.
The readers are assumed to be familiar with the rudiments of geometric
Gromov--Witten theory.

We start with a very brief review of geometric Gromov--Witten theory
in Section~\ref{s:2}, mainly to fix the necessary notations.
In Section~\ref{s:3} the \emph{genus zero axiomatic Gromov--Witten theory} is
introduced. Among the important properties of the axiomatic theory is the 
following theorem:
The ``moduli space'' of the genus zero axiomatic theories of a fixed rank
is acted upon by the twisted loop group.
Furthermore, the subspace of semisimple theories is a homogeneous space of
the twisted loop group.
Then, in Section~\ref{s:4}, the semisimple genus zero theories are 
\emph{quantized} to obtain the higher genus theories.
Here, the above theorem of genus zero theories plays an important role.
The implications of the axiomatic formulation to \emph{Virasoro constraints} 
are discussed in Section~\ref{s:5}.
Finally, the notion of \emph{invariance of tautological equations} and its
applications are briefly discussed in Section~\ref{s:6}.

\emph{Warning}: This article does not aim to give a historical account of 
Gromov--Witten theory, which is preferably left to other experts.
It rather emphasizes upon some highlights centered at the axiomatic 
theory which captures our imagination.
Therefore, the choice of the topics is entirely personal, 
and some important progress is completely left out when its
intersection with axiomatic theory is minimal.

\begin{acknowledgments}
I wish to thank all my collaborators on this subject, D.~Arcara, A.~Givental,
R.~Pandharipande, for many discussions.
This research is partially supported by NSF and an AMS Centennial Fellowship.
\end{acknowledgments}

\section{Review of geometric Gromov--Witten theory} \label{s:2}
Gromov--Witten theory studies the tautological intersection theory on
$\Mbar_{g,n}(X,\beta)$, the moduli stacks of stable maps from curves $C$ of
genus $g$ with $n$ marked points to a smooth projective variety $X$. The
intersection numbers, or \emph{Gromov--Witten invariants}, are integrals of
tautological classes over the virtual fundamental classes of
$\Mbar_{g,n}(X,\beta)$
\[
  \int_{[\Mbar_{g,n}(X,\beta)]^{\vir}}
    \prod_{i=1}^n \ev_i^*(\gamma_i) \psi_i^{k_i}.
\]
Here $\gamma_i \in H^*(X)$
and $\psi_i$ are the \emph{cotangent classes} (gravitational
descendents).

For the sake of the later reference, let us fix some notations.
\begin{enumerate}
\item[(i)] $H := H^*(X, \QQ)$ is a $\QQ$-vector space, assumed of rank $N$.
Let $\{ \phi_{\mu} \}_{\mu=1}^N$ be a basis of $H$.

\item[(ii)] $H$ carries a symmetric bilinear form, Poincar\'e pairing,
\[
  \langle \cdot, \cdot \rangle : H \otimes H \to \QQ.
\]
Define 
\[
  g_{\mu \nu} := \langle \phi_{\mu}, \phi_{\nu} \rangle
\]
and $g^{\mu \nu}$ to be the inverse matrix.

\item[(iii)] Let $\mathcal{H}_t := \oplus_{k=0}^{\infty} H$ be the infinite
dimensional complex vector space with basis $\{ \phi_{\mu} \psi^k \}$.
$\mathcal{H}_t$ has a natural $\QQ$-algebra structure:
\[
 \phi_{\mu} \psi^{k_1} \otimes \phi_{\nu} \psi^{k_2} \mapsto 
 (\phi_{\mu} \cdot \phi_{\nu})  \psi^{k_1 + k_2},
\]
where $\phi_{\mu} \cdot \phi_{\nu}$ is the cup product in $H$.

\item[(iv)] Let $\{ t^{\mu}_k \}$, $\mu=1, \ldots, N$, $k=0, \ldots, \infty$,
be the dual coordinates of the basis $\{ \phi_{\mu} \psi^k \}$.

\end{enumerate}

We note that at each marked point, the insertion is $\mathcal{H}_t$-valued. 
Let
\[
 t:= \sum_{k, \mu} t^{\mu}_k \phi_{\mu} \psi^k
\] 
denote a general element in the vector space $\mathcal{H}_t$.

\begin{enumerate}

\item[(v)] Define 
\[ 
 \la \p^{\mu_1}_{k_1} \ldots \p^{\mu_n}_{k_n}
 \ra_{g,n,\beta} := \int_{[\Mbar_{g,n}(X,\beta)]^{\vir}} \prod_{i=1}^n
 \ev_i^*(\phi_{\mu_i}) \psi_i^{k_i}
\]
and define 
\[
 \la t^n \ra_{g,n,\beta}=\la t \ldots t \ra_{g,n,\beta}
\]
by multi-linearity.

\item[(vi)] Let 
\[
  F^X_g(t) := \sum_{n, \beta} \frac{1}{n!} \la t^n \ra_{g,n,\beta}
\]
be the generating function of all genus $g$ Gromov--Witten invariants.
\footnote{In Gromov--Witten theory, one usually has to deal with 
the coefficients in the \emph{Novikov ring}, due to some convergence issues.
We shall not touch upon this subtleties here but refer the readers to
\cite{LP}.}
The \emph{$\tau$-function of $X$} is the formal expression 
\begin{equation} \label{e:1}
 \tau_{GW}^X := e^{\sum_{g=0}^{\infty} \hbar^{g-1} F_g^X}.
\end{equation}
\end{enumerate}

\section{Genus zero axiomatic Gromov--Witten theory} \label{s:3}

Let $H$ be a $\QQ$-vector space of dimension $N$ with a distinguished element
$\mathbf{1}$.
Let $\{ \phi_\mu \}$ be a basis of $H$ and $\phi_{\mathbf{1}} = \mathbf{1}$. 
Assume that $H$ is endowed with a nondegenerate symmetric $\QQ$-bilinear form,
or metric, $\langle \cdot,\cdot \rangle$. 
Let $\HH$ denote the infinite dimensional vector space $H[z,z^{-1}]$ 
consisting of Laurent polynomials with coefficients in $H$. 
\footnote{Different completions of $\HH$ are used in different places. 
Although there is not a single completion which works for all theorems quoted
in this context, the final results nonetheless make sense as a coherent theory.
This subtlety will be not be discussed in the present article.
See \cite{LP} for the details.} 
Introduce a symplectic form $\Omega$ on $\HH$:
\[
  \Omega (f(z),g(z) ) := \on{Res}_{z=0} \langle f(-z), g(z) \rangle,
\]
where the symbol $\on{Res}_{z=0}$ means to take the residue at $z=0$.

There is a natural polarization $\HH = \HH_{q}\oplus \HH_{p}$ by the 
Lagrangian subspaces $\HH_{q} := H[z]$ and $\HH_{p} := z^{-1} H [ z^{-1} ]$
which provides a symplectic identification of $(\HH, \Omega)$ with the 
cotangent bundle $T^*\HH_{q}$ with the natural symplectic structure.
$\HH_{q}$ has a basis 
\[
 \{ \phi_{\mu} z^k \}, \quad 1 \le \mu \le N, \quad 0 \le k
\]
with dual coordinates $\{ q^k_\mu \}$.
The corresponding basis for $\HH_{p}$ is
\[
 \{ \phi_{\mu} z^{-k-1} \}, \quad 1 \le \mu \le N, \quad 0 \le k
\]
with dual coordinates $\{ p^k_\mu \}$.

For example, let $\{ \phi_{i} \}$ be an \emph{orthonormal} basis of $H$. 
An $H$-valued Laurent formal series can be written in this basis as
\begin{multline*}
 \ldots + (p^1_1,\ldots,p^N_{1}) \frac{1}{(-z)^2}
 + (p^1_{0},\ldots, p^N_{0}) \frac{1}{(-z)} \\
 + (q^1_{0},\ldots, q^N_{0})
 + (q^1_{1}, \ldots, q^N_{1}) z + \ldots.
\end{multline*}
In fact, $\{ p_k^{i}, q_k^{i} \}$ for $k= 0, 1, 2, \ldots$ and
$i=1, \ldots, N$ are the Darboux coordinates compatible with this
polarization in the sense that
\[
 \Omega = \sum_{i,k} d p^{i}_k \wedge d q^{i}_k .
\]

The parallel between $\HH_q$ and $\HH_t$ is evident, and is in fact
given by the following affine coordinate transformation, 
called the \emph{dilaton shift},
\[
  t^{\mu}_k = q^{\mu}_k + \delta^{\mu \mathbf{1}} \delta_{k 1}.
\]

\begin{definition} \label{d:1} 
Let $G_0(t)$ be a (formal) function on $\HH_t$. 
The pair $T:=(\HH, G_0)$ is called a \emph{$g=0$ axiomatic theory} if 
$G_0$ satisfies three sets of genus zero
tautological equations: 
the {\em Dilaton Equation} \eqref{e:de},
the {\em String Equation} \eqref{e:se} and 
the \emph{Topological Recursion Relations} (TRR) \eqref{e:g0trr}.

\begin{align}
 \label{e:de}
 &\frac{\p G_0(t)}{\p t^{\mathbf{1}}_1} (t)
 = \sum_{k=0}^{\infty} \sum_{\mu} t^{\mu}_k
 \frac{\p G_0(t)}{\p t^{\mu}_k} - 2G_0 (t), \\
  \label{e:se}
  &\frac{\p G_0 (t)}{\p t^{\mathbf{1}}_0} =
  \frac{1}{2} \langle t_0,t_0 \rangle + \sum_{k=0}^{\infty}
  \sum_{\nu} t_{k+1}^{\nu} \frac{\p G_0 (t)}{\p t^{\nu}_k},  \\
  \label{e:g0trr}
  &\frac{\p^3 G_0 (t)}
       {\p t^{\alpha}_{k+1} \p t^{\beta}_{l} \p t^{\gamma}_m}
  = \sum_{\mu \nu} \frac{\p^2 G_0 (t)}{\p t^{\alpha}_{k} \p t^{\mu}_{0}}
    g^{\mu \nu} \frac{\p^3 G_0 (t)}
       {\p t^{\nu}_{0} \p t^{\beta}_{l} \p t^{\gamma}_m}, 
  \quad \forall \alpha, \beta, \gamma, k,l,m.
\end{align}
\end{definition}
To simplify the notations, $p_k$ will stand for the vector
$(p^1_{k},\ldots,p^N_{k})$ and $p^{\mu}$ for
$(p^{\mu}_0, p^{\mu}_1, \ldots )$.
Similarly for $q, t$.

In the case of geometric theory, $G_0 = F_0^X$ 
It is well known that $F_0^X$ satisfies the above three
sets of equations \eqref{e:de} \eqref{e:se} \eqref{e:g0trr}.
The main advantage of viewing the genus zero theory through this formulation, 
seems to us, is to replace $\mathcal{H}_t$ by $\HH$ where a symplectic 
structure is available. 
Therefore many properties can be reformulated in terms of the symplectic
structure $\Omega$ and hence independent of the choice of the polarization.
This suggests that the space of genus zero axiomatic Gromov--Witten
theories, i.e. the space of functions $G_0$ satisfying the string equation, 
dilaton equation and TRRs, has a huge symmetry group.

\begin{definition} \label{d:2}
Let $L^{(2)}GL(H)$ denote the {\em twisted loop group} which consists of
$\operatorname{End}(H)$-valued formal Laurent series $M(z)$ in
the indeterminate $z^{-1}$ satisfying $M^*(-z)M(z)=\mathbf{I}$.
Here $\ ^*$ denotes the adjoint with respect to $(\cdot ,\cdot )$.
\end{definition}

The condition $M^*(-z)M(z)=\mathbf{I}$ means that $M(z)$ is a 
symplectic transformation on $\HH$.

\begin{theorem} \label{t:1} \cite{aG4}
The twisted loop group acts on the space of axiomatic genus zero theories.
Furthermore, the action is transitive on the semisimple theories of a fixed
rank $N$.
\end{theorem}

\begin{remarks*}
(i) In the geometric theory, $F^X_0(t)$ is usually a formal function in $t$.
Therefore, the corresponding function in $q$ would be formal at 
$q = - \mathbf{1}z$.
Furthermore, the Novikov rings are usually needed to ensure the 
well-definedness of $F^X_0(t)$. (cf.~Footnote~1.)

(ii) It can be shown that the axiomatic genus zero theory over complex numbers
is equivalent to the definition of abstract (formal) Frobenius manifolds, 
not necessarily conformal. The coordinates on the
corresponding Frobenius manifold is given by the following map \cite{DW}
\begin{equation} \label{e:dg}
 s^{\mu} := \frac{\p}{\p t^{\mu}_0}\frac{\p}{\p t^{\mathbf{1}}_0} G_0 (t).
\end{equation}
\emph{From now on, the term ``genus zero axiomatic theory'' is identified 
with ``Frobenius manifold''}.

(iii) The above formulation (or the Frobenius manifold 
formulation) does not include the divisor axiom, which is true for any
geometric theory.

(iv) Coates--Givental \cite{CG} and Givental \cite{aG4} gave a beautiful 
geometric reformation of the genus axiomatic theory in terms of 
\emph{Lagrangian cones} in $\HH$.
When viewed in the Lagrangian cone formulation, Theorem~\ref{t:1} 
becomes transparent and a proof is almost unnecessary.

Roughly, the descendent Lagrangian cones are constucted in the following way.
Denote by $\mathcal{L}$ the graph of the differential $dG_0$:
\[ 
 \mathcal{L} = \{ ({p}, q) \in T^*\HH_{q}:\ 
 p^{\mu}_k = \frac{\partial}{\partial q^{\mu}_k} G_0 \}.
\]
It is considered as a formal germ at $q = -z$ (i.e.~$t=0$) of a Lagrangian
section of the cotangent bundle $T^*\HH_{q} = \HH$, due to the convergence 
issues of $G_0$. 
$\mathcal{L}$ is therefore considered as a formal germ of a Lagrangian 
submanifold in the space $(\HH,\Omega)$.

\begin{theorem*} 
$(\HH, G_0)$ defines an axiomatic theory if the corresponding
Lagrangian cone $\mathcal{L} \subset \HH$ satisfies the following properties:
$\mathcal{L}$ is a Lagrangian cone with the vertex at the origin of $q$
such that its tangent spaces $L$ are tangent to $\mathcal{L}$ 
exactly along $zL$.
\end{theorem*}

A Lagrangian cone with the above property is also called
\emph{over-ruled} (descendent) Lagrangian cones.
\end{remarks*}

\section{Quantization and higher genus axiomatic theory} \label{s:4}

\subsection{Preliminaries on quantization}

To quantize an infinitesimal symplectic transformation, or its corresponding
quadratic hamiltonians, we recall the standard Weyl quantization. 
A polarization $\HH=T^* \HH_q$ on the symplectic vector space $\HH$ (the phase
space) defines a configuration space $\HH_q$. 
The quantum ``Fock space'' will be a certain class of functions 
$f(\hbar, q)$ on $\HH_q$ (containing at least polynomial functions), 
with additional formal variable $\hbar$ (``Planck's constant''). 
The classical observables are certain functions of $p, q$. The
quantization process is to find for the classical mechanical system on $\HH$ a
``quantum mechanical'' system on the Fock space such that the classical
observables, like the hamiltonians $h(q,p)$ on $\HH$, are quantized to become
operators $\widehat{h}(q,\dfrac{\p}{\p q})$ on the Fock space.

Let $A(z)$ be an $\on{End}(H)$-valued Laurent formal series in $z$ satisfying
\[
  (A(-z) f(-z), g(z)) + (f(-z), A(z) g(z)) =0,
\]
then $A(z)$ defines an infinitesimal symplectic transformation
\[
  \Omega(A f, g) + \Omega(f, A g)=0.
\]
An infinitesimal symplectic transformation $A$ of $\HH$ corresponds to a
quadratic polynomial 
\footnote{Due to the nature of the infinite dimensional vector spaces involved,
the ``polynomials'' here might have infinite many terms, but the degrees 
remain finite.}
$P(A)$ in $p, q$
\[
  P(A)(f) := \frac{1}{2} \Omega(Af, f) .
\]

Choose a Darboux coordinate system $\{ q^i_k, p^i_k \}$. 
The quantization $P \mapsto \widehat{P}$ assigns
\begin{equation} \label{e:wq}
 \begin{split}
  &\widehat{1}= 1, \  \widehat{p_k^{i}}= \sqrt{\hbar} \frac{\p}{\p q^{i}_k}, \
   \widehat{q^{i}_k} = q^{i}_k / {\sqrt{\hbar}}, \\
  &\widehat{p^{i}_k p^{j}_l} = \widehat{p^{i}_k} \widehat{p^{j}_l}
    =\hbar \frac{\p}{\p q^{i}_k} \frac{\p}{\p q^{j}_l}, \\
   &\widehat{p^{i}_k q^{j}_l} = q^{j}_l \frac{\p}{\p q^{i}_k},\\
  &\widehat{q^{i}_k q^{j}_l} = {q}^{i}_k {q}^{j}_l /\hbar ,
 \end{split}
\end{equation}
In summary, the quantization is the process
\[
\begin{matrix}
  A   &\mapsto & P(A)  &\mapsto & \widehat{P(A)} \\
  \text{inf. sympl. transf.}  &\mapsto & \text{quadr. hamilt.}  
    &\mapsto & \text{operator on Fock sp.}.
\end{matrix}
\]
It can be readily checked that the first map is a Lie algebra isomorphism:
The Lie bracket on the left is defined by $[A_1, A_2]=A_1 A_2 - A_2 A_1$ and 
the Lie bracket in the middle is defined by Poisson bracket
\[
 \{ P_1(p,q), P_2(p,q) \} = \sum_{k,i} 
 \frac{\p P_1}{\p p^{i}_k} \frac{\p P_2}{\p q^{i}_k}
 -\frac{\p P_2}{\p p^{i}_k} \frac{\p P_1}{\p q^{i}_k}.
\]
The second map is not a Lie algebra homomorphism, but is very close to being
one.
\begin{lemma} \label{l:1}
\[
 [\widehat{P_1},\widehat{P_2}] =  \widehat{\{ P_1, P_2 \}} + 
  \mathcal{C}(P_1,P_2),
\]
where the cocycle $\mathcal{C}$, in orthonormal coordinates, vanishes except
\[
 \mathcal{C}(p^{i}_k p^{j}_l, q^{i}_k q^{j}_l) = 
 -\mathcal{C}(q^{i}_k q^{j}_l, p^{i}_k p^{j}_l) 
 = 1 + \delta^{i j} \delta_{kl}.
\]
\end{lemma}

\begin{example*}
Let $\dim H=1$ and $A(z)$ be multiplication by $z^{-1}$. 
It is easy to see that $A(z)$ is infinitesimally symplectic.
\begin{equation} \label{e:cse}
 \begin{split}
 P(z^{-1})= &-\frac{q_0^2}{2} - \sum_{m=0}^{\infty} q_{m+1} p_m \\
 \widehat{P(z^{-1})} = &-\frac{q_0^2}{2} 
            - \sum_{m=0}^{\infty} q_{m+1} \frac{\p}{\p q_m}.
 \end{split}
\end{equation}
\end{example*}

Note that one often has to quantize the symplectic instead of the 
infinitesimal symplectic transformations. 
Following the common practice in physics, define
\begin{equation} \label{e:q}
  \widehat{e^{A(z)}} := e^{\widehat{A(z)}},
\end{equation}
for $e^{A(z)}$ an element in the twisted loop group.

\subsection{$\tau$-function for the axiomatic theory}

Let $X$ be the space of $N$ points and $H^{N pt}:= H^*(X)$.
Let $\phi_{i}$ be the delta-function at the $i$-th point.
Then $\{ \phi_{i} \}_{i=1}^N$ form an orthonormal basis 
and are the idempotents of the quantum product
\[
  \phi_{i} * \phi_{j} = \delta_{ij} \phi_{i}.
\]
The genus zero potential for $N$ points is nothing but a sum of genus zero 
potentials of a point
\[
  F^{N pt}_0 (t^1, \ldots, t^N) = F^{pt}_0 (t^1) + \ldots + 
  F^{pt}_0(t^N).
\]
In particular, the genus zero theory of $N$ points is semisimple.

By Theorem~\ref{t:1}, any \emph{semisimple} genus zero axiomatic theory $T$ 
of rank $N$ can be obtained from $H^{N pt}$ by action of an element $O^T$ 
in the twisted loop group. 
By Birkhoff factorization, $O^T = S^T (z^{-1}) R^T (z)$, where
$S(z^{-1})$ (resp.~$R(z)$) is a matrix-valued function in $z^{-1}$
(resp.~$z$). 

In order to define the axiomatic higher genus potentials $G_g^T$ for the 
semisimple theory $T$, one first introduces the ``$\tau$-function of $T$''.

\begin{definition} \label{d:3} \cite{aG2}
Define the \emph{axiomatic $\tau$-function} as
\begin{equation} \label{e:taug}
  \tau_G^T:= \widehat{S^T} (\widehat{R^T} \tau^{N pt}_{GW}),
\end{equation}
where $\tau^{N pt}_{GW}$ is defined in \eqref{e:1}.
Define the \emph{axiomatic genus $g$ potential} $G_g^T$ via the formula 
(cf.~\eqref{e:1})
\begin{equation} \label{e:2}
  \tau_G^T =: e^{\sum_{g=0}^{\infty} \hbar^{g-1} G_g^T}.
\end{equation}
\end{definition}

\begin{remarks*}
(i) It is not obvious that the above definitions make sense.
The function $\widehat{S^T} (\widehat{R^T} \tau^{N pt})$ 
is well-defined, due to some finiteness properties of $\tau^{pt}$, called 
the $(3g-2)$-jet properties \cite{aG2}\cite{eG4}.
The fact that $\log \tau_G^T$ can be written as 
$\sum_{g=0}^{\infty} \hbar^{g-1} (\text{formal function in $t$})$ is also
nontrivial.
The interested readers are referred to the original article \cite{aG2}
or \cite{LP} for details.

(ii) What makes Givental's axiomatic theory especially attractive are the 
facts that
\begin{enumerate}
\item[(a)] It works for any semisimple Frobenius manifolds, not necessarily
coming from geometry.
\item[(b)] It enjoys properties often complementary to the geometric theory.
\end{enumerate}
These will be put into use in the following sections.
\end{remarks*}

\section{Virasoro constraints} \label{s:5}

\subsection{Virasoro operators for points via quantization}
Let $H^{N pt}$ be the genus zero theory for $X$ being $N$ points.
Define the differential operators $D$ on the corresponding 
$\HH =H^{N pt}((z^{-1}))$
\[
  D := z (\frac{d}{dz}) z = z^2 \frac{d}{dz} + z.
\]
Define the operators $\{ L_m \}$ for $m= -1,0,1,2, \ldots$
\[
  L_m := - z^{-1/2} D^{m+1} z^{-1/2}.
\]
The operators $L_m$ have only integer exponents of $z$.

\begin{lemma} \label{l:2}
\begin{itemize}
\item[(i)] 
\begin{equation} \label{e:11}
  [L_m, L_n]=(m-n)L_{m+n}.
\end{equation}
\item[(ii)] $L_m$ are infinitesimal symplectic transformations.
\end{itemize}
\end{lemma}

\begin{proof}
Part (i) can be proved in the following way. 
First perform a change of variables $w = 1/z$. 
Then 
\[
 z^{1/2} L_m z^{-1/2} = (-1)^m  w \frac{d^{m+1}}{d w^{m+1}}.
\]
The RHS has a Fourier transform to the standard vector 
fields on the disk $(-1)^m w^{m+1}\frac{d}{d w}$, 
which obviously satisfies the Virasoro relations of (i).

Part (i) implies that the Lie algebra spanned by $L_m$ is generated by
$L_2$ and $L_{-1}$. These two operators can be verified to satisfy (ii).
\end{proof}

Due to Lemma~\ref{l:2} (ii), $L_m$ can be quantized to $\widehat{L_m}$.
The Lie algebra generated by $\{ \widehat{L_m} \}_{m \ge -1}$ satisfies
the Virasoro relations due to Lemma~\ref{l:1}. 
The Virasoro operators $\{ \widehat{L_m} \}_{m \ge -1}$ constructed above
are the same as $N$ copies of those used in Witten's conjecture in
relation to KdV hierarchies \cite{eW1}.

\subsection{Virasoro operators for semisimple axiomatic theories}
Now for any axiomatic theory of rank $N$, one may define the 
Virasoro operators. The notations in Definition~\ref{d:3} will be followed.

\begin{definition} \label{d:4} \cite{aG2}
\begin{equation} \label{e:12}
 \widehat{L_m^T} := \widehat{S^T} \big( \widehat{R^T} \widehat{L^{H^{N pt}}_m}
   \widehat{R^T}^{-1} \big) \widehat{S^T}^{-1}.
\end{equation}
\end{definition}

\begin{lemma} \label{l:3}
\begin{itemize}
\item[(i)] $\{\widehat{L_m^T}\}_{m\ge -1}$ satisfy the Virasoro relations 
\eqref{e:11}.
\item[(ii)] $\widehat{L_m^T} \tau^T_G = 0$.
\end{itemize}
\end{lemma}

\begin{proof}
(i) is obviously true as the conjugation does not change the commutation 
relations.
(ii) follows from the fact that $\widehat{L^{H^{N pt}}_m} \tau^{H^{N pt}} =0$,
which is $N$ copies of Witten's conjecture.
\end{proof}

\subsection{Virasoro constraints}
\begin{Virasoro} \cite{EHX}
For any projective manifold $X$, there exist ``Virasoro operators''
$\{\widehat{L_m^X}\}_{m\ge -1}$, satisfying the relations \eqref{e:11}, 
such that 
\[
  \widehat{L_m^X} \tau^X_{GW} =0 , \ \forall m \ge -1.
\]
\end{Virasoro}
A good reference for a precise statement can be found in \cite{eG3}.

With Lemma~\ref{l:3}, a clear path to prove Virasoro conjecture would be 
to show, when $T^X$, the genus zero Gromov--Witten theory of $X$, is 
semisimple, the following two statements.
\begin{itemize}
\item[(a)] Definition~\ref{d:4} of Virasoro operators coincide with
the definition of Eguchi--Hori--Xiong \cite{EHX} in the semisimple case.
\item[(b)] $\tau^X_{GW} = \tau^{T^X}_G$.
\end{itemize}

\begin{remarks*}
(i) (a) can be proved with some efforts. See \cite{aG2} and \cite{LP}.
In fact, the operators defined by \eqref{e:12} are also equivalent to
those defined by Dubrovin--Zhang \cite{DZ2}.

(ii) Givental uses the axiomatic framework and a clever observation to
give a one-line proof of genus zero Virasoro Conjecture \cite{aG1}.
\end{remarks*}

(b) will be called \emph{Givental's Conjecture}.
The proof of (b) is more complicated.

Assume that $X$ has a torus action.
Then the equivariant Gromov--Witten invariants are defined and hence
the $\tau$-function $\tau^X_{eGW}$ for the equivariant Gromov--Witten theory
of $X$.
Suppose furthermore that the torus action has isolated fixed points and 
one-dimensional orbits.
It is not hard to see that the genus zero equivariant Gromov--Witten 
invariants define a semisimple axiomatic theory $T^X_e$ \cite{LP}.
Therefore the above formulation works in this case and

\begin{theorem} \cite{aG1, aG2}
Givental's Conjecture holds in the above (equivariant) case.
\end{theorem}

Givental's proof of this theorem is a tour de force and
involves deep theory of Frobenius manifolds and localization.
See \cite{LP} for details.

If the non-equivariant genus zero Gromov--Witten theory $T^X$ is semisimple,
then it is reasonable to expect that non-equivariant limit of $T^X_e$ exists 
and equals $T^X$.
This will prove (non-equivariant) Virasoro Conjecture in a great deal of
examples.
Unfortunately, this proves harder than one expects.
Givental provided a key reduction to a much simpler statement of the existence
of the non-equivariant limit of a certain $R^{T^X_e}$ restricted to the 
``small quantum cohomology'',
which we will not state but refer the interested readers to \cite{aG2}.
This last statement has been proved for the toric Fano manifolds \cite{aG2},
for general toric manifolds \cite{hI}, 
and for some classical flag varieties \cite{JK} \cite{BCK1} \cite{BCK2}.

Another approach to Givental's conjecture is through the invariance of the
tautological equations. See Section~\ref{s:6} for an explanation of the 
following result.

\begin{theorem} \label{t:3} \cite{GL, ypL2, ypL3}
Givental's Conjecture holds for genus one and two.
\end{theorem}

In 2005, C.~Teleman announces a very strong classification theorem of
all semisimple Frobenius manifolds. 
This is a very exciting progress as Givental's Conjecture follows as a 
corollary.

\begin{theorem} \label{t:Teleman} \cite{cT}
Givental's Conjecture is true. Therefore, Virasoro constraints also
holds for semisimple Gromov--Witten theory. 
Furthermore, the equality holds at the level of cycles.
\end{theorem}

\section{Invariance of tautological relations} \label{s:6}

\subsection{Tautological rings}
A basic reference for tautological rings is \cite{rV},
where the history of the subject is explained.

The tautological rings $R^*(\Mbar_{g,n})$ are subrings of
$A^*(\Mbar_{g,n})_{\QQ}$, or subrings of $H^{2*}(\Mbar_{g,n})$ via cycle maps,
generated by some ``geometric classes'' which will be described below.

There are two types of natural morphisms between moduli stacks of curves.
The \emph{forgetful morphisms}
\begin{equation} \label{e:ft}
 \ft_i: \Mbar_{g,n+1} \to \Mbar_{g,n}
\end{equation}
forget one of the $n+1$ marked points.
The \emph{gluing morphisms}
\begin{equation} \label{e:gl}
 \Mbar_{g_1,n_1 +1} \times \Mbar_{g_2, n_2 +1} \to \Mbar_{g_1 + g_2, n_1+n_2},
 \quad \Mbar_{g-1,n+2} \to \Mbar_{g,n},
\end{equation}
glue two marked points to form a curve with a new node.
Note that the boundary strata are the images (of the repeated applications)
of the gluing morphisms, up to factors in $\QQ$ due to automorphisms.

\begin{definition}
The system of tautological rings $\{R^*(\Mbar_{g,n})\}_{g,n}$ is the
smallest system of $\QQ$-unital subalgebra
closed under the forgetful and gluing morphisms.
\end{definition}
As it contains the units, the fundamental classes of the boundary strata
are contained in $R^*(\Mbar_{g,n})$.
From some elementary manipulations, one can also produce Chern classes
of certain tautological vector bundles: $\psi$-classes,
$\lambda$-classes and $\kappa$-classes.

\subsection{Invariance constraints} \label{s:6.2}
Due to the existence of natural \emph{stabilization morphisms}
\begin{equation} \label{e:st}
  \on{st}: \Mbar_{g,n}(X, \beta) \to \Mbar_{g,n}
\end{equation}
defined by forgetting morphisms to $X$ and stabilizing the domain curves 
if necessary, any relations in $H^*(\Mbar_{g,n})_{\QQ}$ can be pull-backed to 
$H^*(\Mbar_{g,n}(X, \beta))_{\QQ}$.
Therefore, the induced equation will hold for any geometric Gromov--Witten 
theory. 
It is natural to expect that the same is true for any axiomatic theory.

For any axiomatic theory $T$, $\tau$-function can be obtained as 
$\tau^T_G = \widehat{S^T} (\widehat{R^T} \tau^{N pt}_{GW})$.
The $\tau$-function of $N$ points certainly satisfies any induced equations.
Therefore, in order to show that the induced tautological relations
hold for $T$, it is enough to show that these relations are \emph{invariant}
under the action of $\widehat{R^T}$ and $\widehat{S^T}$.

\begin{theorem} \label{t:4} \cite{ypL2, ypL3, FSZ}
The tautological relations are invariant under the action of twisted
loop groups.
\end{theorem}

Since the twisted loop groups can be (Birkhoff) factorized into $\widehat{S}$ 
and $\widehat{R}$, one can prove Theorem~\ref{t:4} by proving invariance
individually. 
The $\widehat{S}$ invariance was shown in \cite{ypL2, ypL3} to be a 
consequence of geometry on moduli of stable maps.
The proof of $\widehat{R}$ invariance was originally given in \cite{ypL2},
using Theorem~\ref{t:Teleman}.
\footnote{In fact, only the ``cycle-form'' of Theorem~\ref{t:Teleman} for
$X=\PP^1$ was used. 
The $X=\PP^1$ case can be proved via localization and clever resummation trick
of Givental \cite{aG2}, \emph{without} Teleman's result.
However, a localization proof does not yet exist in the literature.
(It is part of project \cite{LP}.)}
It was later independently discovered by Faber--Shadrin--Zvonkine, and
R.~Pandharipande and the author, that the $\widehat{R}$ invariance
is a simple consequence of geometry on moduli of curves. 
See Section~3 of \cite{FSZ} for details.

\begin{remark*}
The invariance under the action of $\widehat{S}$ imposes little restriction.
It is a consequence of the fact that there is an accountable difference
\cite{KM} \cite{aG2} between $\psi$-classes defined on $\Mbar_{g,n}(X,\beta)$ 
and the pull-backs of $\psi$-classes from  $\Mbar_{g,n}$ via
the stabilization morphism \eqref{e:st}.
On the other hand, the invariance under the action of $\widehat{R}$ imposes
very strong constraints on the structure of tautological relations, and hence
can be used to derive many tautological relations.
\end{remark*}

\subsection{Application I: Virasoro constraints}
As alluded in Theorem~\ref{t:3}, the invariance conjecture can be used
to prove Givental's Conjecture for $g \le 2$.
The idea is to first show some uniqueness theorem: If two genus $g$ potentials
$F_g$ and $G_g$ satisfy enough tautological relations, then
they are identical up to some initial conditions.
For $g=1$ and $2$, this has been done in \cite{DZ2} and \cite{xL1}.
Therefore, Givental's conjecture would follow from the statement that
$G_g$ satisfy those tautological relations.
As explained above, that in turns follows from the invariance of the
tautological relations under the action of $\widehat{R}$.
That is how Theorem~\ref{t:3} is proved.
One immediate corollary of Theorem~\ref{t:3} is Virasoro conjecture in 
$g \le 2$ in the semisimple cases.

\subsection{Application II: Witten's (generalized) conjecture}
Another consequence, which seems less obvious but more interesting to us,
is the following conjecture of Witten.

Witten has proposed a conjectural relation between invariants on the moduli 
space of higher spin curves and the Gelfand--Dickey hierarchies \cite{eW2},
generalizing his previous conjecture \cite{eW1}.
It turns out that the genus zero case defines a semisimple theory,
which is isomorphic to the Frobenius structure of the miniversal
deformation of $A_{r}$-singularities.
The axiomatic theory built on that satisfies the Gelfand--Dickey hierarchies
\cite{aG3}.
Although the ``geometric'' theory defined by the invariants on the moduli 
space of higher spin curves is, strictly speaking, 
not a geometric Gromov--Witten theory,
it also fit into the framework of the axiomatic theory.
Therefore the uniqueness theorems apply and one has the following corollary.

\begin{theorem} \cite{FSZ}
Witten's generalized conjecture is true.
\end{theorem}

This would be a simple consequence of Theorem~\ref{t:Teleman}, combined with
an earlier result by Givental \cite{aG3}.
However, the proof in \cite{FSZ} is independent of Teleman's results.
It uses the ingredients outlined in \cite{ypL0, ypL3}, while expertly 
combining with ingredients from well-known results in tautological rings, 
which we failed to see four years ago.

\subsection{Application III: Finding tautological relations} \label{s:6.5}
Another application is going the opposite direction: a feedback from
Gromov--Witten theory to tautological rings on moduli of curves.
Using Theorem~\ref{t:4}, all the known tautological relations are
obtained by requiring invariance under $\widehat{R}$.
It seems plausible that $\widehat{R}$ invariance can be used to derive
many tautological relations.
In fact, combining Theorem~\ref{t:4} with some known results 
(e.g.~Betti number calculation of Getzler in \cite{eG2}), all 
known tautological relations can be proved uniformly within this framework.
See \cite{ypL1} \cite{AL1} \cite{AL2} \cite{GL} for discussions and 
computations.

\section{Other applications}

\subsection{Twisted Gromov--Witten invariants}
In the axiomatic theory, Coates and Givental \cite{CG} was able to
find the right framework to express the twisted Gromov--Witten invariants
in terms of the untwisted Gromov--Witten invariants.

Given a vector bundle $E$ on $X$, one can define a element $\mathcal{R}$
in the $K$-theory of $\Mbar:=\Mbar_{g,n}(X,\beta)$ in the following way.
Let $\pi: \mathcal{C} \to \Mbar$ be the universal curve and 
$f: \mathcal{C} \to X$ be the universal morphism.
$\mathcal{R}:= R^0 \pi_* f^* E \to R^1 \pi_* f^* E$ is then an element in 
the $K$-theory of $\Mbar$.
The twisted Gromov--Witten invariants are 
the integration, over the virtual fundamental classes $[\Mbar]^{\vir}$
of the usual insertions $H^*(X)$ and descendents, but ``twisted'' with 
certain multiplicative characteristic classes of the virtual vector bundle
$\mathcal{R}$.

It is known in \cite{FP} that these twisted classes can be expressed 
in terms of untwisted classes in the case $E =\mathcal{O}$. 
The underlying geometry was the 
Grothendieck--Riemann--Roch calculation in \cite{dM}.
The case $E$ being any virtual bundle proceed similarly without difficulty.
However, there was no ``closed form'' to express this relation.
The axiomatic framework made this possible.
Of course, it is a truly tour de force to carry out this program even within
the axiomatic framework.

\subsection{Crepant resolution conjecture}
In another direction, this framework is applied to formulate a precise
conjecture regarding the relations between the Gromov--Witten invaraiants
of an orbifold and those of its crepant resolutions.
A recent preprint \cite{tCR} gives an excellent account of this conjecture, 
so we will be very brief.
\footnote{The only \emph{difference} of our points of view lies on the our 
choices of the (over-rulded) Lagrangian cones.
While the authors of \cite{tCR} insists in using descendent cones,
we think the ancestor cones would make the statement a lot cleaner.}

A morphism $f:Y \to \mathcal{X}$ is called a \emph{crepant resolution}, if
$Y$ is smooth and $\mathcal{X}$ is $\QQ$-Gorenstein (e.g. an orbifold) 
such that $f^*(K_{\mathcal{X}}) = K_Y$.
In the case $\mathcal{X}$ is an orbifold, there is a well-defined
orbifold Gromov--Witten theory due to Chen--Ruan \cite{wCR}.
The \emph{crepant resolution conjecture} asserts a close relation between
the Gromov--Witten theory of $Y$ and that of $\mathcal{X}$.
Under suitable conditions (e.g.~hard Lefschetz), 
the conjecture asserts that the Gromov--Witten theory of $Y$ and 
$\mathcal{X}$ are identified, up to analytic continuation on the Novikov
variables.
When the orbifold $\mathcal{X}$ does not satisfy the hard Lefschetz condition,
it is speculated that a very weak relation still hold.
Coates--Corti--Iritani--Tseng and Ruan believe that this weak relation
could and should only be formulated in terms of a symplectic transformation
of the special type.
The interested readers might consult \cite{tCR, CCIT} and references therein, 
as well as the next subsection, for more information.

\subsection{Invariance of Gromov--Witten theory under simple flops}
Crepant resolution is a special case of $K$-equivalence.
$X$ and $X'$ are called $K$-equivalent if there is a common resolution $Y$
such that $K_X$ and $K_{X'}$ are equal after pulling back to $Y$.
In \cite{LLW}, the case when $X$ and $X'$ are smooth and related by a simple 
flop is studied.
It was shown that the big quantum rings are isomorphic \emph{after}
an analytic continuation in the quantum variables.

In \cite{ILLW}, axiomatic framework is used to generalize this result to 
higher genus.
It is first shown, by degeneration to the normal cone, that the statement 
can be reduced to the a statement about toric varieities and toric flops.
Then the quantization formulation reduces the higher genus statement
to a genus zero statement.
An explicit calculation was carried out in genus zero to show the invariance 
of \emph{ancestor} invariants, after analytic continuation.
In fact, results in \cite{ILLW} suggest that the ancestor formulation
is the right framework to study crepant resolution conjecture.

\section{Final remarks}

\subsection{Orbits of twisted loop group action on the moduli spaces of 
Frobenius manifolds}
We learned from Theorem~\ref{t:1} that the twisted loop groups acts on the
space of all Frobenius manifolds of a given rank.
Furthermore, the semisimple theories lie in a single orbit.
It is obvious that the group action has to preserve the 
``degree of diagonizability'' (or degree of semisimplicity) of the theories.
However, it is not known (to me) how many orbits one is to have with a fixed
degree of diagonizability.
It would be very interesting to investigate the orbit structure of the 
twisted loop group action.

\subsection{Integrable hierarchies}
In the proof of the Virasoro conjecture, a basic trick is to move
$\widehat{L_m^{N pt}}$ to the left of $\widehat{S^T} \widehat{R^T}$ and 
to identity it as the $\widehat{L_m^T}$.
In the same spirit, since the $\tau^{N pt}$ satisfies the Hirota equations
of KdV hierarchies, one may try to commute Hirota operators of KdV hierarchies
with $\widehat{S^T} \widehat{R^T}$.
The result could be some yet unknown integrable hierarchies.
This would realize Dubrovin's program of finding a correspondence between
certain classes of integrable hierarchies and the Gromov--Witten theory.
However, the task of commuting operators are rather difficult, due to
some convergence issues.
Successful examples include \cite{aG3} and \cite{GT}.

\subsection{Axiomatic relative/orbifold/open-string theory}
Givental's theory has so far been developed upon the ``original'' 
Gromov--Witten theory.
It is reasonable to ask whether this theory can be generalized to cover
some ramifications of the original Gromov--Witten theory, such as
\emph{relative Gromov--Witten theory} 
or \emph{orbifold Gromov--Witten theory}.
My personal guess is that orbifold theory is probably easier.
Indeed, it is not very difficult to see that Givental's theory of quantization
should work for orbifold theory as well.
At least in the equivariant context, the original localization scheme seems to 
work. 
See Tseng's work \cite{hhT} for some progress along this direction.

In a proper sense that the orbifold compactification is a ``minimal'' one
inside the relative compactification, which involves bubbling off the target
spaces and is more complicated. 
\footnote{This was discovered in a discussion with Y.~Ruan more than three 
years ago.
One example fully worked out is the comparison of the TQFT developed by 
Bryan--Pandharipande in arXiv:math/0411037 (relative compactification) 
and one developed by R.~Cavalieri in arXiv:math/0411500 and arXiv:math/0512225 
(orbifold compactification).
More recently, the relationship between the relative invariants and orbifold
invariants have been extensively studied by D.~Abramovich, C.~Cadman, 
B.~Fantechi, and J.~Wise.}
Therefore, one expects that the axiomatization of relative theory will be
harder.

Furthermore, Gromov--Witten theory is considered as a topological field theory
associated to closed strings. There is an open string analogue.
It is also reasonable to ponder the possibility of an open-string axiomatic
theory, whose geometric theory has not been successfully constructed.
This might further our understanding of the open-string GW theory.

\end{document}